\documentclass[11pt]{jsg}

\usepackage{natbib}
\usepackage{enumerate}
\usepackage{cleveref}
\usepackage{stmaryrd}

        \theoremstyle{plain} 
        \newtheorem{theorem}             {Theorem}  [section]

        \theoremstyle{definition}
        \newtheorem{definition} [theorem]{Definition}

        \theoremstyle{remark}
        \newtheorem{remark}              {Remark}

\begin{document}

\title[Discrete Hamilton--Pontryagin mechanics on Lie groupoids]
{Discrete Hamilton--Pontryagin mechanics and generating functions on
  Lie groupoids}

\author{Ari Stern}

\address{Department of Mathematics\\
University of California, San Diego\\
9500 Gilman Drive \#0112\\
La Jolla, California 92093-0112}

\email{astern@math.ucsd.edu}

\begin{abstract}
  We present a discrete analog of the recently introduced
  Hamilton--Pontryagin variational principle in Lagrangian mechanics.
  This unifies two, previously disparate approaches to discrete
  Lagrangian mechanics: either using the discrete Lagrangian to define
  a finite version of Hamilton's action principle, or treating it as a
  symplectic generating function.  This is demonstrated for a discrete
  Lagrangian defined on an arbitrary Lie groupoid; the often
  encountered special case of the pair groupoid (or Cartesian square)
  is also given as a worked example.
\end{abstract}


\maketitle

\section{Introduction}

In a recent paper, \citet{YoMa2006b} introduced the
Hamilton--Pontryagin variational principle for Lagrangian mechanics.
Given a smooth configuration manifold $Q$ and a Lagrangian $ L \colon
T Q \rightarrow \mathbb{R} $, this principle defines an action for
paths in the so-called Pontryagin bundle $ T Q \oplus T ^\ast Q $,
whose elements are written $ \left( q, v, p \right) $.  Critical paths
for this action, in addition to satisfying the usual Euler--Lagrange
equations $ \dot{p} = \partial L (q,v) / \partial q $, also satisfy
the the Legendre transform $ p = \partial L (q,v) / \partial v $ and
the second-order curve condition $ \dot{q} = v $.  This can be seen as
a unification of two equivalent, but previously disparate, approaches
to Lagrangian mechanics, where the Lagrangian is used either (a) to
define Hamilton's action functional, studying its critical paths; or
(b) to define the Legendre transform, using it to pull back the
canonical symplectic structure from the cotangent bundle to the
tangent bundle.

There is a similar dilemma in the approach to discrete Lagrangian
mechanics.  Given a discrete Lagrangian $ L _h \colon Q \times Q
\rightarrow \mathbb{R} $ (or more generally, $ L _h \colon G
\rightarrow \mathbb{R} $, where $ G \rightrightarrows Q $ is a Lie
groupoid), one can either (a) use $ L _h $ to define a discrete
version of Hamilton's action principle, or (b) treat $ L _h $ as a
symplectic generating function.  While some progress has been made
towards combining these approaches into a discrete
Hamilton--Pontryagin principle for certain important special
cases---namely, the pair groupoid $ Q \times Q $ when $Q$ is either a
vector space \citep{KhYaToKaMaScDe2006,LeOh2008} or a Lie group
\citep{BoMa2009}---this problem has not yet been resolved for general
configuration manifolds $Q$, nor for the even more general case of a
Lie groupoid $ G \rightrightarrows Q $.

In this paper, we present a discrete Hamilton--Pontryagin principle,
which is shown to unify these two disparate approaches (variational
principles vs. generating functions) for discrete Lagrangian mechanics
on arbitrary Lie groupoids.  We begin, in \cref{sec:lagrangian}, by
giving a brief review of Lagrangian mechanics, including the
continuous Hamilton--Pontryagin principle, as well as summarizing the
existing frameworks for discrete Lagrangian mechanics on $ Q \times Q
$ and on Lie groupoids $ G \rightrightarrows Q $.  Next, in
\cref{sec:discreteHP}, we introduce the discrete Hamilton--Pontryagin
principle, which is defined with respect to paths in the cotangent
groupoid $ T ^\ast G \rightrightarrows A ^\ast G $ beginning at the
zero section; the approach is related to that used by
\citet{Milinkovic1999} in studying the Morse homology of generating
functions.  This variational principle and its solutions, which imply
those of the previous approaches to discrete Lagrangian mechanics, are
derived first in cotangent bundle coordinates and then given
intrinsically.  Finally, in \cref{sec:pairGroupoid}, we work out the
special case of the pair groupoid $ Q \times Q $, in particular
showing how this corresponds to the formulation of \citet{LeOh2008}
when $Q$ is a vector space.

\section{Lagrangian Mechanics, Continuous and Discrete}
\label{sec:lagrangian}

\subsection{Lagrangian and Hamilton--Pontryagin mechanics}
Let $Q$ be a smooth configuration manifold, and $ L \colon T Q
\rightarrow \mathbb{R} $ be a Lagrangian on its tangent bundle.  There
are two theoretically equivalent, but conceptually distinct,
approaches to the Lagrangian mechanics of this system.  The first,
which we will call the variational approach, is to study critical
paths $ q \colon [a,b] \rightarrow Q $ for the
action functional
\begin{equation*}
  S (q) = \int _a ^b L \left( q (t) , \dot{q} (t)
  \right) \mathrm{d}t .
\end{equation*}
Such a path is critical if and only if it satisfies the
Euler--Lagrange equations
\begin{equation*}
  \frac{\mathrm{d}}{\mathrm{d}t} \frac{ \partial L }{ \partial \dot{q}
  } =   \frac{ \partial L }{ \partial q }  .
\end{equation*}
The second approach, which we will call the symplectic approach, is to
use the Legendre transform $ \mathbb{F} L \colon T Q \rightarrow
T^\ast Q $, which takes $ (q,v) \mapsto (q, \partial L / \partial v)
$, to pull back the canonical symplectic structure from $ T ^\ast Q $
to $ T Q $.  If $ \omega $ is the canonical symplectic $2$-form on $ T
^\ast Q $, then one can define a $2$-form $ \omega _L = \left(
  \mathbb{F} L \right) ^\ast \omega $ on $ T Q $.  (The form $ \omega
_L $ is symplectic when $L$ is a hyperregular Lagrangian, and
presymplectic more generally.)  Next, one defines the energy function
$ E \colon T Q \rightarrow \mathbb{R} $ given by
\begin{equation*}
  E \left( q, v \right) = \mathbb{F}  L \left( q, v \right) \cdot v -
  L \left( q , v \right) .
\end{equation*}
Finally, then, one looks for vector fields $X \in \mathfrak{X} \left(
  T Q \right) $ that satisfy
\begin{equation*}
  i _X \omega _L = \mathrm{d} E ,
\end{equation*}
which is essentially the tangent bundle version of Hamilton's
equations on $ T ^\ast Q $.  (For further background, see
\citealp[Chapter 7]{MaRa1999}.)

\citet{YoMa2006b} showed that these two approaches can be unified
through an expanded variational principle, which they call the {\em
  Hamilton--Pontryagin principle}.  Given a path $ \left( q, v, p
\right) \colon [a,b] \rightarrow T Q \oplus T ^\ast Q $, the
Hamilton--Pontryagin action is given by
\begin{equation*}
  \tilde{ S } \left( q, v, p \right) = \int _a ^b  \left[ L \left( q (t), v (t)
    \right) + p (t) \cdot \left( \dot{q} (t) - v (t) \right) \right]
  \mathrm{d}t .
\end{equation*}
This is essentially the usual action functional---except, rather than
simply prescribing the second-order curve constraint $ \dot{q} = v $,
one treats $q$ and $v$ as independent variables and then uses $p$ as a
Lagrange multiplier to enforce this constraint.  Varying over paths
with prescribed endpoints, so that $ \delta q (a) = 0 $ and $ \delta q
(b) = 0 $, the variation of the action is
\begin{multline*}
  \mathrm{d} \tilde{ S } \left( q, v, p \right) \cdot \left( \delta q
    , \delta v , \delta p \right) \\
  \begin{aligned}
    &= \int _a ^b \left[ \frac{ \partial L }{ \partial q } \cdot \delta
      q + \frac{ \partial L }{ \partial v } \cdot \delta v + p \cdot
      \left( \delta \dot{q} - \delta v \right) + \delta
      p \cdot \left( \dot{q} - v \right) \right] \mathrm{d} t \\
    &= \int _a ^b \left[ \left( \frac{ \partial L }{ \partial q } -
        \dot{p} \right) \cdot \delta q + \left( \frac{ \partial L
        }{ \partial v } - p \right) \cdot \delta v + \delta p \cdot
      \left( \dot{q} - v \right) \right] \mathrm{d} t .
  \end{aligned}
\end{multline*} 
Therefore, $ \left( q, v, p \right) $ is a critical path if and only
if it solves the so-called {\em implicit Euler--Lagrange equations}
\begin{equation*}
  \dot{p} = \frac{ \partial L }{ \partial q } , \qquad p =
  \frac{ \partial L }{ \partial v } , \qquad \dot{q} = v ,
\end{equation*}
which combine the Euler--Lagrange equations, the Legendre transform,
and the second-order curve condition into a single set of equations.
It should be noted that this principle is especially useful for
studying constrained and other degenerate systems, and is closely
connected with the generalized Legendre transform of
\citet{Tulczyjew1974,Tulczyjew1977}, the generalized Hamiltonian
dynamics formalism of \citet{SkRu1983}, and the Dirac structures of
\citet{Courant1990}.

\subsection{Discrete Lagrangian mechanics} The idea of discrete
Lagrangian mechanics was put forward in seminal papers by
\citet{Suris1990} and \citet{MoVe1991}, among others, and a general
theory was developed over the subsequent decade.  (See
\citealp{MaWe2001}, for a comprehensive survey.)  This work was
motivated by the need to develop structure-preserving (e.g.,
symplectic) numerical integrators for Lagrangian mechanical systems on
general configuration manifolds; the methods developed using this
discrete Lagrangian framework are called {\em variational
  integrators}.

For most of the work in this field, the starting point is to replace
the Lagrangian $ L \colon T Q \rightarrow \mathbb{R} $ by a {\em
  discrete Lagrangian} $ L _h \colon Q \times Q \rightarrow \mathbb{R}
$, which approximates the contribution to the action integral,
\begin{equation*}
  L _h \left( q _0, q _1 \right) \approx \int _{ t _0 } ^{ t _1 } L
  \left( q (t) , \dot{q} (t) \right) \mathrm{d}t ,
\end{equation*}
for a time step of size $ h = t _1 - t _0 $.  As with continuous
Lagrangian mechanics, there are two typical ways to proceed, following
either the variational or the symplectic point of view.

The variational approach to discrete Lagrangian mechanics is as
follows.  Suppose that we specify a sequence of time steps $ a = t _0
< t _1 < \cdots < t _N = b $, with equal step size $ h = t _{ n + 1 }
- t _n $ for $ n = 0, \ldots, N - 1 $.  We then define a discrete
path to be a sequence of configuration points $ q _0, \ldots , q
_N \in Q$; this can be thought of as approximating a continuous path $
q \colon [a,b] \rightarrow Q $, with $ q _n \approx q \left( t _n
\right) $.  Given the discrete Lagrangian $ L _h \colon Q \times Q
\rightarrow \mathbb{R} $, the {\em discrete action sum} is defined to
be
\begin{equation*}
  S _h \left( q _0, \ldots, q _N \right) = \sum _{ n = 0 } ^{ N
    - 1 } L _h \left( q _n , q _{ n + 1 } \right)  \approx \int _a ^b
  L \left( q (t) , \dot{q} (t) \right) \mathrm{d}t .
\end{equation*}
Next, taking fixed-endpoint variations of the discrete path, so that $
\delta q _0 = 0 $ and $ \delta q _N = 0 $, it follows that
\begin{multline*}
  \mathrm{d} S _h \left( q _0 , \ldots, q _N \right) \cdot \left(
    \delta q _0, \ldots , \delta q _N \right) \\
  = \sum _{ n = 1 } ^{ N - 1 } \left[ \partial _0 L _h \left( q _n , q
      _{ n + 1 } \right) + \partial _1 L _h \left( q _{ n - 1 } , q _n
    \right) \right] \cdot \delta q _n .
\end{multline*}
Therefore, a discrete path is critical if and only if it satisfies the
{\em discrete Euler--Lagrange equations}
\begin{equation*}
  \partial _0 L _h \left( q _n , q _{ n + 1 } \right) + \partial _1 L
  _h \left( q _{ n - 1 } , q _n \right) = 0 , \qquad n = 1, \ldots
  , N - 1 .
\end{equation*}
This implicitly defines a two-step numerical integrator on $ Q \times
Q $, which (given suitable assumptions of nondegeneracy) maps $ \left(
  q _{ n - 1 } , q _n \right) \mapsto \left( q _n , q _{ n + 1 }
\right) $.

On the other hand, the symplectic approach to discrete Lagrangian
mechanics is to view $ L _h \colon Q \times Q \rightarrow \mathbb{R} $
as the generating function for a symplectic map on $ T ^\ast Q $.  To
do this, one defines the {\em discrete Legendre transforms} $
\mathbb{F} ^{ \pm } L _h \colon Q \times Q \rightarrow T ^\ast Q $ by
\begin{equation*}
  \mathbb{F} ^- L _h \left( q _0, q _1 \right) = - \partial _0 L _h \left(
    q _0, q _1 \right) , \qquad   \mathbb{F} ^+ L _h \left( q _0, q _1
  \right) =  \partial _1 L _h \left( q _0, q _1 \right).
\end{equation*}
Therefore, we can implicitly define a map on $ T ^\ast Q $ by
\begin{equation*}
  p _0 =   \mathbb{F}  ^- L _h \left( q _0, q _1 \right), \qquad p _1 =
  \mathbb{F} ^+  L _h \left( q _0, q _1 \right) ,
\end{equation*}
and if $ \mathbb{F} ^- L _h $ is invertible, this defines one step of
the symplectic integrator
\begin{equation*}
  \mathbb{F}  ^+ L _h \circ \left( \mathbb{F} ^- L _h \right) ^{-1} \colon
  T^\ast Q \rightarrow T^\ast Q , \qquad  \left( q _0, p _0 \right)
  \mapsto \left( q _1 , p _1 \right) .
\end{equation*}
More precisely, the discrete Legendre transforms define the Lagrangian
submanifold of $ \left( T ^\ast Q , - \omega \right) \times \left( T
  ^\ast Q , \omega \right) $ generated by $ L _h $. From this
perspective, if the invertibility condition holds, then this
submanifold is the graph of a symplectic map on $T^\ast Q$
\citep[see][]{Weinstein1971,Weinstein1979}.

Note that, if we perform the composition in the opposite order, we get
the previously derived two-step method,
\begin{equation*}
  \left( \mathbb{F} ^- L _h \right) ^{-1} \circ \mathbb{F} ^+ L _h \colon Q
  \times Q \rightarrow Q \times Q, \qquad \left( q _{ n - 1 } , q _n
  \right) \mapsto \left( q _n , q _{ n + 1 } \right) ,
\end{equation*} 
where the discrete Euler--Lagrange equations follow automatically from
the fact that $ \mathbb{F} ^- L _h \left( q _n , q _{ n + 1 } \right)
= \mathbb{F} ^+ L _h \left( q _{ n - 1 }, q _n \right) $.

\subsection{Discrete Lagrangian mechanics and Lie groupoids}
\citet{Weinstein1996} observed that both approaches in the previous
section can be generalized using Lie groupoids.  Let $ G
\rightrightarrows Q $ be a given Lie groupoid, and define a discrete
Lagrangian $ L _h \colon G \rightarrow \mathbb{R} $.  The earlier
formulations then coincide with the special case $ G = Q \times Q $,
which is called the {\em pair groupoid}.  This perspective has
continued to bear fruit in recent years, being further developed by
\citet{MaMaMa2006} and extended to discrete nonholonomic Lagrangian
mechanics by \citet{IgMaMaMa2008}.

In the variational approach, one begins by taking a fixed element $ g
\in G $ and considering the space of {\em admissible sequences}, which
consist of composable elements $ g _1, \ldots, g _N \in G $ such that
$ g _1 \cdots g _N = g $.  The discrete action for an admissible
sequence is then taken to be
\begin{equation*}
  S _h \left( g _1, \ldots, g _N \right) = \sum _{ n = 1 } ^N L _h
  \left( g _n \right) ,
\end{equation*}
and discrete Euler--Lagrange equations are obtained by finding
sequences which are critical for this action function.  In the case of
the pair groupoid $ Q \times Q $, fixing an element of the groupoid
corresponds simply to fixing the endpoints $ g = \left( q _0, q _N
\right) $, while the set of admissible sequences $ \left( q _0, q _1
\right) , \ldots, \left( q _{ N-1}, q _N \right) \in Q \times Q $ can
be identified with the sequence of configuration points $ q _0,
\ldots , q _N \in Q $.

\citet{MaMaMa2006} showed that the discrete Euler--Lagrange equations
can be expessed in terms of left- and right-invariant vector fields on
$G$, each of which can be identified with sections of the Lie
algebroid $ A G \rightarrow Q $ associated to $G$.  To describe this,
we must first introduce some notation.  Let $ \alpha, \beta \colon G
\rightarrow Q $ denote the source and target maps, $ \epsilon \colon Q
\rightarrow G $ denote the identity section, and $ i \colon G
\rightarrow G $ denote the inversion map on $G$.  For any $ g \in G $,
define the left- and right-translation maps, respectively, by
\begin{align*}
  \ell _g \colon \alpha ^{-1} \left( \beta (g) \right) &\rightarrow
  \alpha ^{-1} \left( \alpha (g) \right) & r _g \colon \beta ^{-1}
  \left( \alpha (g) \right) &\rightarrow \beta ^{-1}
  \left( \beta (g) \right)  \\
  g ^\prime &\mapsto g g ^\prime & g ^\prime &\mapsto g ^\prime g .
\end{align*}
Then, given a section $ X \in \Gamma \left( A G \right) $, the
left-invariant vector field $ \overleftarrow{ X } \in \mathfrak{X} (G)
$ is defined by $ \overleftarrow{ X } (g) = T _{ \epsilon \left( \beta
    (g) \right) } \ell _g \left( X \left( \beta (g) \right) \right) $.
Similarly, the right-invariant vector field $ \overrightarrow{ X } \in
\mathfrak{X} (G) $ is given by $ \overrightarrow{ X } (g) = - \left( T
  _{ \epsilon \left( \alpha (g) \right) } \left( r _g \circ i \right)
\right) \left( X \left( \alpha (g) \right) \right) $.  This implies
the following relationship between the brackets $ \llbracket \cdot,
\cdot \rrbracket $ on $ \Gamma \left( A G \right) $ and $ [ \cdot,
\cdot ] $ on $ \mathfrak{X} (G) $:
\begin{equation*}
  \overleftarrow{ \llbracket X, Y \rrbracket } = \bigl[ \overleftarrow{ X } ,
    \overleftarrow{ Y } \bigr], \qquad    \overrightarrow{ \llbracket
    X, Y \rrbracket } = - \bigl[  \overrightarrow{ X },
    \overrightarrow{ Y } \bigr] .
\end{equation*}
With these definitions, the discrete Euler--Lagrange equations are
\begin{equation*}
  \overleftarrow{ X } \left[ L _h \right]  \left( g _n \right) =
  \overrightarrow{ X } \left[ L _h \right] \left( g _{ n + 1 } \right)
  , \qquad  n = 1, \ldots, N - 1 ,
\end{equation*}
for all sections $ X \in \Gamma \left( A G \right) $.  In the pair
groupoid case, one sees that
\begin{align*} 
  \overleftarrow{ X } \left[ L _h \right] \left( q _{ n - 1 }, q _n
  \right) &= \partial _1 L _h \left( q _{ n - 1 }
    , q _n \right) \cdot X \left( q _n \right)  \\
  \overrightarrow{ X } \left[ L _h \right] \left( q _n , q _{n+1}
  \right) &= - \partial _0 L _h \left( q _n , q _{ n + 1 } \right)
  \cdot X \left( q _n \right) ,
\end{align*} 
for any $ X \in \mathfrak{X} (Q) $, so this formulation agrees with
the earlier expression of the discrete Euler--Lagrange equations on $
Q \times Q $.

For the symplectic approach, we begin by defining the {\em cotangent
  groupoid} $ T ^\ast G \rightrightarrows A ^\ast G $, where the base
$ A ^\ast G \rightarrow Q $ is the dual of the Lie algebroid $ A G
\rightarrow Q $.  This is a {\em symplectic groupoid}, with a
canonical symplectic structure on $ T ^\ast G $ and Poisson structure
on $ A ^\ast G $; moreover, the source and target maps are
anti-Poisson and Poisson, respectively \citep{CoDaWe1987,Marle2005}.
Explicitly, the source map $ \tilde{ \alpha } \colon T ^\ast G
\rightarrow A ^\ast G $ and target map $ \tilde{ \beta } \colon T
^\ast G \rightarrow A ^\ast G $ may be defined such that, taking $ \mu
\in T _g ^\ast G $,
\begin{equation*}
  \tilde{ \alpha } \left( \mu \right) \cdot X \left( \alpha (g)
  \right) = \mu \cdot \overrightarrow{ X } (g) , \qquad \tilde{ \beta
  } \left( \mu \right) \cdot X \left( \beta (g) \right) = \mu \cdot
  \overleftarrow{ X } (g) ,
\end{equation*} 
for every section $ X \in \Gamma \left( A G \right) $.
(The multiplicative structure of $ T ^\ast G $ will not be necessary
here, so we omit discussion of it.)

Now, the discrete Lagrangian $ L _h \colon G \rightarrow \mathbb{R} $
generates a Lagrangian submanifold of the cotangent groupoid, $
\mathrm{d} L _h (G) \subset T ^\ast G $.  If $ L _h $ is suitably
nondegenerate, it follows that this submanifold determines a Poisson
automorphism $ A ^\ast G \rightarrow A ^\ast G $
\citep{Weinstein1996,MaMaMa2006}.  Specifically, define the discrete
Legendre transforms $ \mathbb{F} ^\pm L _h \colon G \rightarrow A
^\ast G $ by
\begin{equation*}
  \mathbb{F}  ^- L _h = \tilde{ \alpha } \circ \mathrm{d} L _h ,
  \qquad \mathbb{F}  ^+ L _h = \tilde{ \beta } \circ \mathrm{d} L _h .
\end{equation*}
Therefore, $ \left\{ \left( \mathbb{F} ^- L _h (g), \mathbb{F} ^+ L _h
    (g) \right) \;\middle\vert\; g \in G \right\} $ is a coisotropic
relation on $ A ^\ast G $, and if $ \mathbb{F} ^- L _h $ is a
diffeomorphism, then this relation is the graph of the discrete flow
map
\begin{equation*}
  \mathbb{F}  ^+ L _h \circ \left( \mathbb{F}  ^- L _h \right) ^{-1}
  \colon A ^\ast G \rightarrow A ^\ast G .
\end{equation*}
To see that this is compatible with the discrete Euler--Lagrange
equations, observe that for any $ X \in \Gamma \left( A G \right) $,
\begin{alignat*}{2}
  \bigl( \tilde{ \alpha } \circ \mathrm{d} L _h \bigr) (g) \cdot X
  \left( \alpha (g) \right) &= \mathrm{d} L _h (g) \cdot
  \overrightarrow{ X } (g) &&= \overrightarrow{ X } \left[ L _h
  \right]
  (g) \\
  \bigl( \tilde{ \beta } \circ \mathrm{d} L _h \bigr) (g) \cdot X
  \left( \beta (g) \right) &= \mathrm{d} L _h (g) \cdot
  \overleftarrow{ X } (g) &&= \overleftarrow{ X } \left[ L _h \right]
  (g) .
\end{alignat*}
Thus, the discrete Euler--Lagrange equations correspond to the
condition $ \mathbb{F} ^+ L _h \left( g _n \right) = \mathbb{F} ^- L
_h \left( g _{ n + 1 } \right) $.  Equivalently, this means that $ g
_1, \ldots, g _N $ is a solution of the discrete Euler--Lagrange
equations when $ \mathrm{d} L _h \left( g _1 \right) , \ldots,
\mathrm{d} L _h \left( g _N \right) $ is a composable sequence in $ T
^\ast G $.

\section{The Discrete Hamilton--Pontryagin Principle}
\label{sec:discreteHP}

For the continuous Hamilton--Pontryagin principle, the key idea was to
relax the condition for curves in $ T Q $ to be second-order, but to
enforce this condition weakly using Lagrange multipliers in $ T ^\ast
Q $.  Analogously, we will develop a discrete Hamilton--Pontryagin
principle by relaxing the requirement that the sequence $ g _1 ,
\ldots, g _N \in G $ is admissible---requiring only that it be {\em
  admissible up to homotopy}---and will weakly enforce the condition
that this homotopy is, in fact, constant.  To allow for the ``Lagrange
multipliers,'' these homotopies will be given by paths in $ T ^\ast G
$ rather than in $G$.
  
To simplify the exposition, we will first sketch this formulation, in
\cref{sec:cotangent}, using cotangent bundle coordinates $ \left( g,
  \mu \right) \in T ^\ast G $.  The fully intrinsic treatment on $ T
^\ast G $ will be given subsequently, in \cref{sec:intrinsic}, along
with definitions and the proof of the main theorem.

\subsection{Formulation in cotangent bundle coordinates}
\label{sec:cotangent}
Suppose that $ \left( g _n , \mu _n \right) \colon [ 0, 1 ]
\rightarrow T ^\ast G $, for $ n = 1, \ldots, N $, is a sequence of
paths in the cotangent groupoid.  Furthermore, require that the
initial point $ \left( g _n (0) , \mu _n (0) \right) $ of each path
lie in the zero section, so that $ \mu _n (0) = 0 _{ g _n (0) } $, and
that $ g _1 (0) , \ldots, g _N (0) $ be an admissible sequence for a
fixed element $ g \in G $.

Given such a sequence of paths, take the discrete Hamilton--Pontryagin
action to be
\begin{equation*}
  \tilde{ S } _h \left( g _1, \mu _1, \ldots, g _N , \mu _N \right) =
  \sum _{ n = 1 } ^N \int _0 ^1 \left[ L _h   \left( g _n (s) \right)
    + \mu _n (s) \cdot g _n ^\prime  (s) \right] \mathrm{d} s .
\end{equation*}
Variations of this action are
\begin{multline*}
  \mathrm{d} \tilde{ S } _h \left( g _1, \mu _1, \ldots, g _N , \mu _N
  \right) \cdot \left( \delta g _1, \delta \mu _1, \ldots, \delta g _N
    , \delta \mu _N \right) \\
  \begin{aligned}
    &= \sum _{ n = 1 } ^N \int _0 ^1 \left[ \mathrm{d} L _h \left( g
        _n (s) \right) \cdot \delta g _n (s) + \mu _n (s) \cdot \delta
      g _n ^\prime (s) + \delta \mu _n (s)
      \cdot g _n ^\prime (s) \right] \mathrm{d} s \\
    &= \sum _{ n = 1 } ^N \biggl( \int _0 ^1 \left[ \left( \mathrm{d}
        L _h \left( g _n (s) \right) - \mu _n ^\prime (s) \right)
      \cdot \delta g _n (s) + \delta \mu _n (s) \cdot g _n ^\prime
      (s) \right] \mathrm{d} s \\
    &\qquad\quad + \mu _n (1) \cdot \delta g _n (1) \biggr).
  \end{aligned}
\end{multline*}
Now, the integral terms vanish when $ \left( g _n , \mu _n \right) $
is a solution of
\begin{equation*}
  \mu _n ^\prime (s) = \mathrm{d} L _h
  \left( g _n (s) \right) , \qquad   g _n ^\prime (s) = 0 .
\end{equation*}
The second equation states that $ g _n (s) = g _n $ is constant, so $
\mu _n (s) $ moves vertically along the fiber $ T _{ g _n } ^\ast G $,
beginning at $ \mu _n (0) = 0 _{ g _n } $.  We use this fact to solve
the remaining equation,
\begin{equation*}
  \mu _n (s) = \mu _n (0) + s \, \mathrm{d} L _h \left( g _n \right) =
  s \, \mathrm{d} L _h \left( g _n \right) ,
\end{equation*}
so in particular, $ \mu _n (1) = \mathrm{d} L _h \left( g _n \right)
$.  Additionally, $ g _1 , \ldots, g _N $ is an admissible sequence,
since $ g _n = g _n (0) $, which was assumed to be admissible.
Finally, restricting to $ \left( g _n , \mu _n \right) $ where these
equations are satisfied, the variation of the action is given by the
remaining boundary terms
\begin{equation*}
  \sum _{ n = 1 } ^N \mu _n (1) \cdot \delta g _n (1) = \sum
  _{ n = 1 } ^N \mathrm{d} L _h \left( g _n \right) \cdot \delta
  g _n .
\end{equation*}
However, the restricted variations $ \delta g _n $ are no longer
completely arbitrary, since they must be tangent to the space of
admissible sequences $ g _1 , \ldots, g _N $.  Therefore, we are back
in the case considered by \citet{MaMaMa2006}, so these terms vanish
precisely when $ g _n $ is a solution of the discrete Euler--Lagrange
equations---or equivalently, when $ \mu _1 (1) , \ldots, \mu _N (1)
\in T ^\ast G $ is a composable sequence.

In summary, $ \mathrm{d} \tilde{ S } _h = 0 $ implies that $ g _1 (1),
\ldots, g _N (1) $ is an admissible sequence in $G$, that $ \mu _n (1)
= \mathrm{d} L _h \left( g _n (1) \right) $, and that $ \mu _1 (1) ,
\ldots, \mu _N (1) $ is a composable sequence in $ T ^\ast G $.

\subsection{Intrinsic formulation}
\label{sec:intrinsic}
We will now make precise the approach sketched out in the previous
section, doing this intrinsically on $ T ^\ast G $ rather than in
cotangent bundle coordinates.  Let $ \pi _G \colon T ^\ast G
\rightarrow G $ be the cotangent bundle projection, $ \tilde{ \theta }
$ be the canonical $1$-form on $ T ^\ast G $, and $ \tilde{ \omega } =
- \mathrm{d} \tilde{ \theta } $ be the canonical symplectic $2$-form.

Given a Lie groupoid $ G \rightrightarrows Q $, define the fibration $
E \rightarrow G $ consisting of paths in $ T ^\ast G $ beginning at
the zero section,
\begin{equation*}
  E = \left\{ \gamma \colon [0,1] \rightarrow T ^\ast G
   \;\middle\vert\; \gamma (0) \in 0 _G \right\} ,
\end{equation*}
where the projection is given by $ \gamma \mapsto \left( \pi _G \circ
  \gamma \right) (1) $.

\begin{definition}
  Given a fixed $ g \in G $, a sequence of paths $ \gamma _1 , \ldots,
  \gamma _N \in E $ is said to be {\em admissible} if the initial
  elements $ \left( \pi _G \circ \gamma _1 \right) (0) , \ldots,
  \left( \pi _G \circ \gamma _N \right) (0) $ form an admissible
  sequence in $G$.
\end{definition}

This implies that the projected sequence $ \left( \pi _G \circ \gamma
  _1 \right) (1) , \ldots , \left( \pi _G \circ \gamma _N \right) (1)
\in G $ is only admissible up to homotopy, which is weaker than the
usual assumption that $ g _1, \ldots, g _N $ must actually be an
admissible sequence.

\begin{definition}
  Let $ L _h \colon G \rightarrow \mathbb{R} $ be a discrete
  Lagrangian, and $ \gamma _1 , \ldots , \gamma _N \in E $ be an
  admissible sequence of paths.  Then the {\em discrete
    Hamilton--Pontryagin action} of this sequence is
\begin{equation*}
  \tilde{ S } _h \left( \gamma _1, \ldots, \gamma _N \right) = \sum _{ n
    = 1 } ^N \left[ \int _0 ^1 \left( L _h \circ \pi _G \right)
    \left( \gamma _n (s) \right) \mathrm{d}s + \int _{ \gamma _n }
    \tilde{ \theta } \right] ,
\end{equation*}
The sequence is said to satisfy the {\em discrete Hamilton--Pontragin
  principle} if
\begin{equation*}
  \mathrm{d} \tilde{ S } _h \left( \gamma _1, \ldots, \gamma _N
  \right) = 0 .
\end{equation*}
\end{definition}

\begin{theorem}
  Let $ \gamma _1 , \ldots, \gamma _N \in E $ be an admissible
  sequence of paths, and denote $ \mu _n = \gamma _n (1) \in T ^\ast G
  $ and $ g _n = \pi _G \left( \mu _n \right) \in G $ for $ n = 1 ,
  \ldots, N $.  If the sequence satisfies the Hamilton--Pontryagin
  principle, then
\begin{enumerate}[(i)]
\item $ g _1, \ldots, g _N $ is an admissible sequence in $G$,
\item $ \mu _n = \mathrm{d} L _h \left( g _n \right) $ for $ n = 1,
  \ldots, N $,
\item $ \mu _1 , \ldots, \mu _N $ is a composable sequence in $ T
  ^\ast G $.
\end{enumerate}
\end{theorem}

These properties may be thought of as, respectively, (i) the discrete
second-order curve condition, (ii) the discrete Legendre transform,
and (iii) the discrete Euler--Lagrange equations.

\begin{proof}
  Take a variation of the Hamilton--Pontryagin action to obtain
\begin{multline*} 
  \mathrm{d} \tilde{ S } _h \left( \gamma _1, \ldots, \gamma _N
  \right) \cdot \left( \delta \gamma _1, \ldots, \delta \gamma _N
  \right) \\
  = \sum _{ n = 1 } ^N \left[ \int _0 ^1 \left( \pi _G ^\ast
      \mathrm{d} L _h \right) \cdot \delta \gamma _n (s) \,\mathrm{d} s
    + \int _{ \gamma _n } \mathfrak{L} _{ \delta \gamma _n } \tilde{
      \theta } \right].
\end{multline*} 
Applying Cartan's ``magic formula,'' the Lie derivative terms become
\begin{align*}
  \int _{ \gamma _n } \mathfrak{L} _{ \delta \gamma _n } \tilde{ \theta }
  &= \int _{ \gamma _n } \left( i _{ \delta \gamma _n} \mathrm{d}
    \tilde{ \theta } + \mathrm{d} i _{ \delta \gamma _n } \tilde{
      \theta } \right) \\
  &= - \int _{ \gamma _n } i _{ \delta \gamma _n } \tilde{ \omega } +
  \int _{ \partial \gamma _n } \tilde{ \theta } \cdot \delta \gamma _n
  \\
  &= \int _0 ^1 i _{ \gamma _n ^\prime (s) } \tilde{ \omega } \cdot
  \delta \gamma _n (s) \,\mathrm{d} s + \tilde{ \theta } \cdot \delta
  \gamma _n (1) .
\end{align*}
Therefore, the variation of the action is
\begin{multline*}
  \mathrm{d} \tilde{ S } _h \left( \gamma _1, \ldots, \gamma _N
  \right) \cdot \left( \delta \gamma _1, \ldots, \delta \gamma _N
  \right) \\
  = \sum _{ n = 1 } ^N \left[ \int _0 ^1 \left( \pi _G ^\ast
      \mathrm{d} L _h + i _{ \gamma _n ^\prime (s) } \tilde{ \omega }
    \right) \cdot \delta \gamma _n (s) \,\mathrm{d}s + \tilde{ \theta
    } \cdot \delta \gamma _n (1) \right].
\end{multline*}
For variations taken along the fibers of $E$, this is stationary when
\begin{equation*}
  i _{ \gamma _n ^\prime (s) } \tilde{
    \omega } = - \pi _G ^\ast \mathrm{d} L _h  , \qquad n = 1, \ldots, N ,
\end{equation*}
i.e., the paths $ \gamma _n $ are solutions to Hamilton's equations
for the singular Hamiltonian $ \tilde{ H } = - L _h \circ \pi _G
\colon T ^\ast G \rightarrow \mathbb{R} $.  However, since $ \tilde{ H
}$ is constant on fibers of $ T ^\ast G $ (which form a Lagrangian
foliation), the Hamiltonian flow is fiber-preserving, and hence the
projection $ \left( \pi _G \circ \gamma _n \right) (s) = g _n $ is
constant.  This proves (i), since $ g _n = \left( \pi _G \circ \gamma
  _n \right) (0) $, which was assumed to be an admissible sequence in
$G$.  Furthermore, the Hamiltonian vector field is constant on fibers
of $ T ^\ast G $, and it follows that $ \gamma _n (s) = s \,\mathrm{d}
L _h \left( g _n \right) $.  Therefore, $ \mu _n = \gamma _n (1) =
\mathrm{d} L _h \left( g _n \right) $, which proves (ii).

Finally, restricting the action to these solutions, the integral terms
vanish, and the remaining boundary terms are
\begin{equation*} 
  \sum _{ n = 1 } ^N \tilde{ \theta } \cdot \delta \gamma _n
  (1) = \sum _{ n = 1 } ^N \mu _n \cdot \delta g _n = \sum _{
    n = 1 } ^N \mathrm{d} L _h \left( g _n \right) \cdot \delta g
  _n .
\end{equation*} 
This is precisely the variation of the usual discrete action, so these
terms vanish when the discrete Euler--Lagrange equations are
satisfied.  Therefore, $ \tilde{ \beta } \left( \mu _n \right) =
\tilde{ \alpha } \left( \mu _{ n + 1 } \right) $ for $ n = 1, \ldots,
N - 1 $, so $ \mu _1, \ldots, \mu _N $ is a composable sequence in $ T
^\ast G $, which completes the proof of (iii).
\end{proof}

\begin{remark}
  For a single time step $ N = 1 $, the action
\begin{equation*}
  \tilde{ S } _h \colon E \rightarrow \mathbb{R} , \qquad \gamma
  \mapsto   \int _0 ^1 \left( L _h \circ \pi _G \right)
  \left( \gamma (s) \right) \mathrm{d}s + \int _{ \gamma }
  \tilde{ \theta }
\end{equation*}
can be understood as a {\em Morse family}, which generates the usual
Lagrangian submanifold $ \mathrm{d} L _h (G) \subset T ^\ast G $.
This is closely related to the work of \citet{Milinkovic1999}, who
studied a similar action principle on $ T ^\ast Q $, in connection
with the Morse homology of generating functions for Lagrangian
submanifolds.  In \citeauthor{Milinkovic1999}'s formulation, the
Lagrangian submanifold is determined by the time-$1$ Hamiltonian
isotopy of $ \tilde{ H } $, which in this case, takes the zero section
to $ \mathrm{d} L _h (G) $.
\end{remark}

\section{Example: The Pair Groupoid}
\label{sec:pairGroupoid}

Let $ G = Q \times Q $ be the pair groupoid, and $ L _h \colon Q
\times Q \rightarrow \mathbb{R} $ be a discrete Lagrangian.  Fixing
the endpoints $ g = \left( q _0, q _N \right) $, suppose that $ \gamma
_1, \ldots, \gamma _N \in E $ is an admissible sequence of paths,
where $ \gamma _n (s) = \left( q _n ^0 (s), p _n ^0 (s) , q _n ^1 (s)
  , p _n ^1 (s) \right) $.  Since the sequence is admissible, this
implies in particular that $ q _0 = q _1 ^0 (0) $, $ q _n ^1 (0) = q
_{ n + 1 } ^0 (0) $ for $ n = 1, \ldots, N - 1 $, and $ q _N ^1 (0) =
q _N $.  The discrete Hamilton--Pontryagin action of this sequence is
then
\begin{multline*}
  \tilde{ S } _h \left( \gamma _1 , \ldots, \gamma _N   \right) \\
  = \sum _{ n = 1 } ^N \int _0 ^1 \left[ L _h \left( q _n ^0
        (s) , q _n ^1 (s) \right) - p _n ^0 (s) \cdot q _n
      ^{0\,\prime} (s) + p _n ^1 (s) \cdot q _n ^{ 1\, ^\prime } (s)
    \right] \mathrm{d} s .
\end{multline*}
Taking variations of this action gives
\begin{multline*}
  \mathrm{d} \tilde{ S } _h \left( \gamma _1 , \ldots, \gamma _N
  \right) \cdot \left( \delta \gamma _1 , \ldots, \delta \gamma _N
  \right) \\
  \begin{aligned}
    &= \sum _{ n = 1 } ^N \int _0 ^1 \Bigl[ \partial _0 L _h \left( q
      _n ^0 (s) , q _n ^1 (s) \right) \cdot \delta q _n ^0 (s)
    + \partial _1 L _h \left( q _n ^0 (s) , q _n ^1 (s) \right)
    \cdot \delta q _n ^1 (s) \\
    & \qquad\qquad\qquad - p _n ^0 (s) \cdot \delta q _n ^{0\,\prime}
    (s) + p _n ^1 (s) \cdot \delta q _n ^{ 1\, ^\prime
    } (s) \\
    & \qquad\qquad\qquad\qquad - \delta p _n ^0 (s) \cdot q _n
    ^{0\,\prime} (s) + \delta p _n ^1 (s) \cdot q _n ^{ 1\, ^\prime }
    (s) \Bigr] \mathrm{d} s ,
    \end{aligned}
\end{multline*}
and integrating by parts, this simplifies to
\begin{multline*}
  \mathrm{d} \tilde{ S } _h \left( \gamma _1 , \ldots, \gamma _N
  \right) \cdot \left( \delta \gamma _1 , \ldots, \delta \gamma _N
  \right) \\
  \begin{aligned}
    &= \sum _{ n = 1 } ^N \biggl( \int _0 ^1 \Bigl[ \left( \partial _0
      L _h \left( q _n ^0 (s) , q _n ^1 (s) \right) + p _n
      ^{0\,\prime} (s) \right) \cdot \delta q _n ^0 (s) \\
    &\qquad\qquad\qquad+ \left( \partial _1 L _h \left( q _n ^0 (s) ,
        q _n ^1 (s) \right) - p _n ^{1\,\prime} (s) \right)
    \cdot \delta q _n ^1 (s) \\
    & \qquad\qquad\qquad - \delta p _n ^0 (s) \cdot q _n
    ^{0\,\prime} (s) + \delta p _n ^1 (s) \cdot q _n ^{ 1\, ^\prime
    } (s) \Bigr] \mathrm{d} s \\
    & \qquad\qquad - p _n ^0 (1) \cdot \delta q _n
    ^0 (1) + p _n ^1 (1) \cdot \delta q _n ^1 (1) \biggr).
    \end{aligned}
\end{multline*}
For variations taken along the fibers of $E$, this is stationary when
the integral terms vanish, so 
\begin{alignat*}{2}
  p _n ^{0\,\prime} (s) &= - \partial _0 L _h \left( q _n ^0 (s) , q
    _n ^1 (s) \right), &\qquad  q _n ^{0\,\prime} (s) &= 0 ,\\
  p _n ^{1\,\prime} (s) &= \partial _1 L _h \left( q _n ^0 (s) , q _n
    ^1 (s) \right), & q _n ^{1\,\prime} (s) &= 0 ,
\end{alignat*}
for $ n = 1 , \ldots, N $.  Therefore, $ q ^0 _n (s) = q ^0 _n $ and $
q ^1 _n (s) = q ^1 _n $ are constant, so 
\begin{equation*}
  p _n ^0 (1) = - \partial _0 L _h \left( q _n ^0 , q _n ^1 \right) ,
  \qquad p _n ^1 (1) = \partial _1 L _h \left( q _n ^0 , q _n ^1
  \right) , \qquad n = 1, \ldots, N .
\end{equation*}
Furthermore, the admissibility assumption implies that $ q _0 = q ^0
_1 $, $ q _n ^1 = q _{ n + 1 } ^0 $ for $ n = 1 , \ldots, N - 1 $, and
$ q _N ^1 = q _N $.

Finally, restricting to these solutions, the restricted variations
must then satisfy $ \delta q _1 ^0 = 0 $, $ \delta q _n ^1 = \delta q
_{ n + 1 } ^0 = $ for $ n = 1, \ldots, N - 1 $, and $ \delta q _N ^1 =
0 $.  Therefore, the remaining terms of the action are
\begin{equation*}
  \sum _{ n = 1 } ^N \left( - p _n ^0 (1) \cdot \delta q ^0 _n + p _n
    ^1 (1) \cdot \delta q _n ^1 \right) =   \sum _{ n = 1 } ^{N-1}
  \left( - p _{ n + 1 } ^0 (1) + p _n ^1 (1) \right) \cdot \delta q _n
  ^1 ,
\end{equation*}
which vanish when $ p _n ^1 (1) = p _{ n + 1 } ^0 (1) $ for $ n = 1,
\ldots, N - 1 $.

In summary, if $ \gamma _1, \ldots, \gamma _N $ satisfies the discrete
Hamilton--Pontryagin principle on $ Q \times Q $, then the following
is true:
\begin{enumerate}[(i)]
\item $ q _0 = q _1 ^0 $, $ q _n ^1 = q _{ n + 1 } ^0 $ for $ n = 1
  ,\ldots, N - 1 $, and $ q _N ^1 = q _N $,
\item $ p _n ^0 (1) = - \partial _0 L _h \left( q _n ^0 , q _n ^1
  \right) $ and $ p _n ^1 (1) = \partial _1 L _h \left( q _n ^0 , q _n
    ^1 \right) $ for $ n = 1, \ldots, N $,
\item $ p _n ^1 (1) = p _{n+1} ^0 (1) $ for $ n = 1 ,\ldots, N - 1 $.
\end{enumerate}
These are, respectively, the discrete second-order curve condition,
the discrete Legendre transform, and the discrete Euler--Lagrange
equations for systems on $ Q \times Q $.

\begin{remark}
  The cotangent paths $ p _n ^0 $, $ p _n ^1 $ can be seen as Lagrange
  multipliers, which serve to enforce the composability conditions $ q
  _n ^1 = q _{n+1} ^0 $, as well as the endpoint constraints $ q _1 ^0
  = q _0 $ and $ q _N ^1 = q _N $.  \citet{LeOh2008} showed that, when
  $Q$ is a vector space (or given a particular local coordinate
  neighborhood), this can be done even more directly, by simply
  subtracting these points rather than taking a path between them.
  This leads to an alternative choice of the action sum,
\begin{equation*}
  \sum _{ n = 1 } ^N L _h \left( q _n ^0 , q _n ^1 \right) + p _0
    \cdot \left( q _1 ^0 - q _0 \right) + \sum _{ n
      = 1 } ^{ N - 1 } p _n \cdot \left( q _{ n + 1 } ^0 - q _n ^1
    \right) + p _N \cdot \left( q _N ^1 - q _N \right) ,
\end{equation*}
whose variations are computed to be
\begin{multline*}
  \sum _{ n = 1 } ^N \left[ \left( \partial _0 L _h \left( q _n ^0 , q
      _n ^1 \right) + p _{n-1} \right) \cdot \delta q _n ^0 +
  \left( \partial _1 L _h \left( q _n ^0 , q _n ^1 \right) - p _n
  \right) \cdot \delta q _n ^1 \right] \\+ \delta p _0 \cdot \left( q _1 ^0 -
    q _0 \right) + \sum _{ n = 1 } ^{ N - 1 } \delta p _n \cdot \left(
    q _{ n + 1 } ^0 - q _n ^1 \right) + \delta p _N \cdot \left( q _N
    ^1 - q _N \right) .
\end{multline*}
Therefore, this vanishes when
\begin{equation*}
  p _{ n - 1 } = - \partial _0 L _h \left( q _n ^0 , q _n ^1 \right) ,
  \qquad p _n = \partial _1 L _h \left( q _n ^0 , q _n ^1
  \right) , \qquad n = 1, \ldots, N 
\end{equation*}
and when 
\begin{equation*}
  q _1 ^0 = q _0 , \qquad   q _N ^1 = q _N , \qquad  q _n ^1 = q _{ n
    + 1 } ^0 , \qquad  n =
  1, \ldots, N - 1 ,
\end{equation*}
which are consistent with the results obtained in this section.

Effectively, the Leok--Ohsawa approach can be interpreted as using a
smaller fibration $ E ^\prime = Q \times Q \times T ^\ast Q \times T
^\ast Q $, consisting only of the endpoints of the paths in $E$, to
define a Morse family over $ Q \times Q $.  However, when $Q$ is not a
vector space, there is no intrinsic, global meaning to subtracting two
points in $Q$.  This action may still be defined locally, in a
neighborhood of the diagonal of $ Q \times Q $, given a choice of
local coordinates or a retraction.  In general, though, one requires
the larger space of paths $E$ in order to define the discrete
Hamilton--Pontryagin principle globally.
\end{remark}

\section*{Acknowledgments}

Thanks to Melvin Leok, Jerry Marsden, David Mart\'in de~Diego, and
Joris Vankerschaver for helpful discussions and feedback during the
development of this work, as well as to the editor and anonymous
referee for their valuable comments and suggestions.  Research
supported in part by NSF (PFC Award 0822283 and DMS Award 0715146), as
well as by NIH, HHMI, CTBP, and NBCR.


\begin{thebibliography}{20}
\providecommand{\natexlab}[1]{#1}

\bibitem[{Bou-Rabee and Marsden(2009)}]{BoMa2009}
Bou-Rabee, N., and J.~E. Marsden (2009), {H}amilton--{P}ontryagin integrators
  on {L}ie groups part {I}: Introduction and structure-preserving properties.
  \emph{Found. Comput. Math.}, \textbf{9} (2), 197--219.
\newblock \href {http://dx.doi.org/10.1007/s10208-008-9030-4}
  {\texttt{doi:10.1007/s10208-008-9030-4}}.

\bibitem[{Coste et~al.(1987)Coste, Dazord, and Weinstein}]{CoDaWe1987}
Coste, A., P.~Dazord, and A.~Weinstein (1987), Groupo\"\i des symplectiques. In
  \emph{Publications du {D}\'epartement de {M}ath\'ematiques. {N}ouvelle
  {S}\'erie. {A}, {V}ol.\ 2}, volume~87 of \emph{Publ. D\'ep. Math. Nouvelle
  S\'er. A}, pages i--ii, 1--62. Univ. Claude-Bernard, Lyon.

\bibitem[{Courant(1990)}]{Courant1990}
Courant, T.~J. (1990), Dirac manifolds. \emph{Trans. Amer. Math. Soc.},
  \textbf{319} (2), 631--661.

\bibitem[{Iglesias et~al.(2008)Iglesias, Marrero, Mart{\'{\i}}n~de Diego, and
  Mart{\'{\i}}nez}]{IgMaMaMa2008}
Iglesias, D., J.~C. Marrero, D.~Mart{\'{\i}}n~de Diego, and E.~Mart{\'{\i}}nez
  (2008), Discrete nonholonomic {L}agrangian systems on {L}ie groupoids.
  \emph{J. Nonlinear Sci.}, \textbf{18} (3), 221--276.

\bibitem[{Kharevych et~al.(2006)Kharevych, Yang, Tong, Kanso, Marsden,
  Schr\"{o}der, and Desbrun}]{KhYaToKaMaScDe2006}
Kharevych, L., W.~Yang, Y.~Tong, E.~Kanso, J.~E. Marsden, P.~Schr\"{o}der, and
  M.~Desbrun (2006), Geometric, variational integrators for computer animation.
  In \emph{SCA '06: Proceedings of the 2006 ACM SIGGRAPH/Eurographics symposium
  on Computer animation}, pages 43--51. Eurographics Association,
  Aire-la-Ville, Switzerland.

\bibitem[{Leok and Ohsawa(2008)}]{LeOh2008}
Leok, M., and T.~Ohsawa (2008), Discrete {D}irac structures and variational
  discrete {D}irac mechanics. Preprint.
\newblock \href {http://arxiv.org/abs/0810.0740}
  {\texttt{arXiv:0810.0740~\![math.SG]}}.

\bibitem[{Marle(2005)}]{Marle2005}
Marle, C.-M. (2005), From momentum maps and dual pairs to symplectic and
  {P}oisson groupoids. In \emph{The breadth of symplectic and {P}oisson
  geometry}, volume 232 of \emph{Progr. Math.}, pages 493--523. Birkh\"auser
  Boston, Boston, MA.

\bibitem[{Marrero et~al.(2006)Marrero, Mart{\'{\i}}n~de Diego, and
  Mart{\'{\i}}nez}]{MaMaMa2006}
Marrero, J.~C., D.~Mart{\'{\i}}n~de Diego, and E.~Mart{\'{\i}}nez (2006),
  Discrete {L}agrangian and {H}amiltonian mechanics on {L}ie groupoids.
  \emph{Nonlinearity}, \textbf{19} (6), 1313--1348.

\bibitem[{Marsden and Ratiu(1999)}]{MaRa1999}
Marsden, J.~E., and T.~S. Ratiu (1999), \emph{Introduction to mechanics and
  symmetry}, volume~17 of \emph{Texts in Applied Mathematics}. Springer-Verlag,
  New York, second edition. A basic exposition of classical mechanical systems.

\bibitem[{Marsden and West(2001)}]{MaWe2001}
Marsden, J.~E., and M.~West (2001), Discrete mechanics and variational
  integrators. \emph{Acta Numer.}, \textbf{10}, 357--514.

\bibitem[{Milinkovi{\'c}(1999)}]{Milinkovic1999}
Milinkovi{\'c}, D. (1999), Morse homology for generating functions of
  {L}agrangian submanifolds. \emph{Trans. Amer. Math. Soc.}, \textbf{351} (10),
  3953--3974.

\bibitem[{Moser and Veselov(1991)}]{MoVe1991}
Moser, J., and A.~P. Veselov (1991), Discrete versions of some classical
  integrable systems and factorization of matrix polynomials. \emph{Comm. Math.
  Phys.}, \textbf{139} (2), 217--243.

\bibitem[{Skinner and Rusk(1983)}]{SkRu1983}
Skinner, R., and R.~Rusk (1983), Generalized {H}amiltonian dynamics. {I}.
  {F}ormulation on {$T^{\ast} Q\oplus TQ$}. \emph{J. Math. Phys.}, \textbf{24}
  (11), 2589--2594.
\newblock \href {http://dx.doi.org/10.1063/1.525654}
  {\texttt{doi:10.1063/1.525654}}.

\bibitem[{Suris(1990)}]{Suris1990}
Suris, Y.~B. (1990), Hamiltonian methods of {R}unge-{K}utta type and their
  variational interpretation. \emph{Mat. Model.}, \textbf{2} (4), 78--87.

\bibitem[{Tulczyjew(1974)}]{Tulczyjew1974}
Tulczyjew, W.~M. (1974), Hamiltonian systems, {L}agrangian systems and the
  {L}egendre transformation. In \emph{Symposia {M}athematica, {V}ol. {XIV}
  ({C}onvegno di {G}eometria {S}implettica e {F}isica {M}atematica, {INDAM},
  {R}ome, 1973)}, pages 247--258. Academic Press, London.

\bibitem[{Tulczyjew(1977)}]{Tulczyjew1977}
Tulczyjew, W.~M. (1977), The {L}egendre transformation. \emph{Ann. Inst. H.
  Poincar\'e Sect. A (N.S.)}, \textbf{27} (1), 101--114.

\bibitem[{Weinstein(1971)}]{Weinstein1971}
Weinstein, A. (1971), Symplectic manifolds and their {L}agrangian submanifolds.
  \emph{Advances in Math.}, \textbf{6}, 329--346.

\bibitem[{Weinstein(1979)}]{Weinstein1979}
Weinstein, A. (1979), \emph{Lectures on symplectic manifolds}, volume~29 of
  \emph{CBMS Regional Conference Series in Mathematics}. American Mathematical
  Society, Providence, R.I. Corrected reprint.

\bibitem[{Weinstein(1996)}]{Weinstein1996}
Weinstein, A. (1996), Lagrangian mechanics and groupoids. In \emph{Mechanics
  day ({W}aterloo, {ON}, 1992)}, volume~7 of \emph{Fields Inst. Commun.}, pages
  207--231. American Mathematical Society, Providence, RI.

\bibitem[{Yoshimura and Marsden(2006)}]{YoMa2006b}
Yoshimura, H., and J.~E. Marsden (2006), Dirac structures in {L}agrangian
  mechanics. {II}. {V}ariational structures. \emph{J. Geom. Phys.}, \textbf{57}
  (1), 209--250.

\end{thebibliography}

\end{document}